\documentclass[12pt,a4paper]{article}
 \usepackage{amsmath,amstext,amssymb,amscd,emlines2 }

\begin{document}
\title{Sequential choice functions and stability problems}
\author{Danilov V.I. \thanks{Central Institute of Economics and Mathematics of the RAS. Nahimovskii prospect 47, 117418, Moscow; email: vdanilov43@mail.ru}}
\date{\today}
\maketitle

The concept of sequential choice functions is introduced and studied. This concept applies to the reduction of the problem of stable matchings with sequential workers to a situation where the workers are linear.

\section{Introduction}

Our interest is related to the concept of stability of contract systems. We assume everywhere, that agent preferences are described by Plott choice functions (Section 2). We additionally assume that the choice functions of some agents are sequential. Sequentiality here means that the choice is made by sequential application of $q$ linear orders or criteria: the first item is selected from the list of available contracts according to the first criterion, the next item is selected from the remaining menu according to the second criterion, and so on. This concept was introduced relatively recently (\cite{CY}, see also \cite{HM}). We show that it is possible to modify the problem in such a way that the solutions do not actually change (so that both problems are equivalent), but agents with sequential choice functions become those with linear choice functions. It is clear that the second problem (with `linear' agents) is easier to analyze (for example, for a construction of rotations).

In the next Section, we recall the definition of Plott choice functions and
introduce the operation $*$ of the sequential application of choice functions. Here some statements about this operation are established, which will be useful in the future. Section 3 recalls the concept of the stability of a matching or a contract system. In Section 4, a design of the modification of the problem is proposed and the main theorems on the reduction of workers are formulated. After some discussion in Section 5, we provide the proofs of the main theorems.

\section{Two constructions of Plott choice functions}

Here we briefly recall the concept of  Plott choice functions and then give two necessary for the further ways to build such functions.

\subsection{Plott choice functions}

Let $X$ be a set. A \emph{choice function} (hereinafter CF) on $X$ is a mapping $C:2^X\to 2^X$, which satisfies the condition $C(A)\subseteq A$ for any $A\subseteq X$. The subset $A$ is called \emph{menu}, and $C(A)$ is called \emph{choice} from this menu. However, we will not be interested in any CFs, but only in those that somehow model the human choice of the best from the available. There is an extensive literature on this topic (see \cite{AA}). As undoubted feature of rational behavior the following property should be recognized, which is most often called \emph{consistency} (there are other names: independence of rejected alternatives, etc.).\medskip

\emph{Consistency}. If $C(A)\subseteq B\subseteq A$, then $C(B)=C(A)$. \medskip

This axiom is usually supplemented by other conditions. For us, the following one will be important, which is called \emph{substitution} (or the heredity axiom): \medskip

\emph{Substitution.} If $B\subseteq A$ then $C(A)\cap B\subseteq C(B)$. \medskip

In other words, if an element from $B$ is selected from the larger set $A$, then it is also selected from the smaller one. In a sense, the items look like as competitors, and appearance of a new item in the menu reduces the chances of other items being selected. This property is not as undoubted as the consistency, but it is very relevant to the stability problems.

We say that a CF is \emph{Plott} one if it is consistent and substitute.
Plott in \cite{P} introduced a class of choice functions, which he called "path independent". These functions were distinguished by the condition
\begin{equation}\label{Plo}
                                 C(A\cup B)=C(C(A)\cup B).                 \end{equation}
Later it was realized that they are exactly the CFs having the properties of
consistency and substitution. The basic properties of such CFs can be found in \cite{AA, M, DK}. Therefore, we will be brief and limit ourselves only to those concepts that will be needed for the future.

A CF $C$ is \emph{cardinally monotonous} if $A\subseteq B$ implies $|C(A)|\le
|C(B)|$.

A menu $A$ is called \emph{acceptable} if $C(A)=A$.

A Plott function $C$ is called \emph{non-empty-valued} if $C(A)\ne \emptyset$ when $A$ is non-empty. For this, it is sufficient to require that all singletons be acceptable.

With every Plott CF $C$ one can associate the \emph{Blair relation} (or the relation of revealed preference) $\preceq =\preceq _C$ on the set $2^X$ by setting
\begin{equation}\label{Bla}
      A\preceq B \text{ if } C(A\cup B)\subseteq B.
\end{equation}
This relation is reflexive and transitive. In the class of menus equivalent to $A$, there is the smallest, namely $C(A)$. 
On the set of acceptable menus, the relation $\preceq$ is antisymmetric, that is a partial order. Moreover, this poset is a lattice in which the join (sup) of two acceptable menus $A$ and $B$ is equal to $C(A\cup B)$.

An element $x$ is called \emph{desirable} (upon the condition $A$) if $x\in C(A\cup x)$. In other words, the agent would like to add $x$ to the menu. In another way, we can say that $A\prec_C A\cup x$, that is, that the menu $A\cup x$ is strictly better (according to Blair relation) than $A$.

We can not do without examples. Here and further, we shall assume that $X$ is finite set. The simplest example of Plott function gives the maximization of a linear order. Let $\le$ be a linear order on the set $X$. Then in each (non-empty) subset $A\subseteq X$ there is a (unique) maximal element $a=\max_\le(A)$. The mapping $A \mapsto  \max_\le(A)$ gives a CF which, as it is easy to understand, is Plott function. We will call such CFs \emph{linear}.

A more interesting and widely used in the literature on stability example of Plott function gives the choice (using a linear order $\le$) of not one best, but the first $q$ best elements, the choice with a \emph{quota} $q$. In \cite{CY} this choice is called \emph{responsive}. It is also a Plott CF. The sequential CF, which will be discussed below, is again the choice of the $q$ best, but with respect to different linear orders. But it is better to start with the more general operation $*$.

\subsection{Sequential choice}

Having several CFs, one can use them to build a new CF. The most popular operation is that of the union. Below we will introduce two more new operations. The main for further is the operation of \emph{sequential (or serial) application} of two (or more) CFs. When we first chose somewhat (and TAKE it away) from the menu using the first CF, and then chose from the rest using the second one. Formally, CF $F*G$ is defined as:
\begin{equation}\label{seq}
                              (F*G)(A)=F(A)\cup G(A-F(A)).
\end{equation}
This operation, as it is easy to understand, is associative, ($F*(G*H)=(F*G) *H$), and therefore, it makes a sense to write $F*G*H*...$ without parentheses. In general case, this the operation is not commutative. \medskip

\textbf{Proposition 1.} \emph{If $F$ and $G$ are Plott CFs, then $F*G$ is a Plott CF.}\medskip

Proof. We need to check the consistency and substitutability of $F*G$.

\emph{Consistency.} By virtue of the definition of $*$, we have $F(A)\cup
G(A-F(A))\subseteq B$. In particular, $F(A)\subseteq B$, from where (using the
consistency of $F$) $F(B)=F(A)$. Next, $G(A-F(B))=G(A-F(A_)$ is also contained in $B$ and does not overlap with $F(B)$. Therefore, $G(A-F(B))\subseteq
B-F(B)\subseteq A-F(B)$, and from the consistency of $G$ we obtain $G(B-F)B))=G(A-F(B))=G(A-F(A))$.

\emph{Substitute.} Let $A\subseteq B$, and let $a\in A$ be chosen from $B$ (by $F*G$); we have to show that $a$ is also chosen from $A$. If $a\in F(B)$ then, from the substitutability of $F$, $a\in F(A)$ and everything is fine. Therefore, we can assume that $a\notin F(A)$ (and even more so it does not belong to $F(B)$) and therefore belongs to $G(B-F(B))$. In addition, $a$ belongs to $A-F(A)$, which is contained in $B-F(B)$. Therefore, from the substitutability of $G$, we get $a\in G(A-F(A))$. $\Box$\medskip

\textbf{Definition.} A CF $C$ is called \emph{sequential} if it is represented as $F_1*...*F_q$, where $F_1,...,F_q$ are linear CFs.\medskip

According to Proposition 1, any sequential CF is a Plott CF. In addition, it is $q$-\emph{quotable}, that is, the choice from $A$ has $q$ elements if $|A|\ge q$, and coincides with $A$ if $|A|\le q$. Any Plott CF with quota 1 is linear. If $F$ and $G$ are quotable with quotas $q$ and $r$, then $F*G$ is quotable with quota $q+r$.

Sequential CFs provide a lot of examples of quotable CFs. However, not all as it shown in \cite{CY}.

\subsection{Integration of Plott functions}

For what follows we need another construction. It can be conditionally called `integration'. This name is given by analogy with the operation of measure integration. In a particular case, this operation has occurred in \cite{DK}.

Suppose we have a mapping $\pi :Y\to X$ of sets. We look at it as a fiber bundle with the base $X$ and the fibers $Y_x=\pi ^{-1}(x)$, where $x$ runs through $X$. In other words, $Y$ is the (disjoint) union of fibers, $Y=\coprod_{x\in X} Y_x$.

Let us assume that, for each $x\in X$, a CF $C_x$ is given on the fiber $Y_x$, and, in addition, a CF $C$ is given on the base $X$. Then one can define CF $D$ by $Y$ (it could be denoted as $\int C_xdC$) using the following formula (where $B\subseteq Y$ and $B_x=B\cap Y_x$):
          \begin{equation}\label{int}
D(B)=\bigcup _{x\in C(\pi (B))} C_x(B_x).
          \end{equation}
In other words, $D(B)$ consists of those $b\in B$ that are selected in fiber $B_x$, if $x=\pi(b)$ is selected in the set $\pi(B)$. For clarity, imagine that there is a multi-storey building $Y$ with several entrances $X$. And there are a set $B$ of available apartments in this house. You would like to occupy some of them. You are primarily interested in entrances, and you first select some (available) entrances. And only then you choose interesting apartments in each of selected entrances.\medskip

\textbf{Example.} Suppose that all CFs ($C_x$ for $x\in X$ and $C$) are linear (with linear orders $\le_x$ and $\le$ correspondingly). Then the `integrated' CF $D$ is linear and corresponds to the ordinal (or lexicographic) sum of these linear orders. \medskip

Next, we will assume that the fiber's CFs $C_x$ are non-empty-valued (in fact, for the further purposes, it is enough to consider linear CFs $C_x$). Under this assumption, the following simple formula is true:
    \begin{equation}\label{pro}
                                 \pi (D(B))=C(\pi (B)).
    \end{equation}
Indeed, the inclusion $\subseteq$ is obvious. On the other hand, if       $x\in C(\pi (B))$, then $B_x$ is not empty, which means that the set $C_x(B_x)$ is also non-empty and is mapped to $x$.

Let us consider a more special case when all fiber's CFs $C_x$ are linear. In this case, the set $C_x(B_x)$ consists of exactly one element (if $B_x\ne \emptyset$), so that $D(B)$ is injectively (and hence bijectively) mapped to its image $C(\pi(B))$. In particular, the size of $D(B)$ is equal to the size of $C(\pi(B))$. Therefore, if the CF $C$ is cardinally monotonous then $D$ is cardinally monotonous as well. For the further, we need the next simple \medskip

\textbf{Lemma 1.} \emph{If all CFs $C_x$ (where $x\in X)$ are linear then the following statements are equivalent:}

           1) $D(B)=B$,

2) \emph{$B$ is bijective to $\pi(B)$, and $\pi(B)=C(\pi(B))$.}  $\Box$ \medskip

\textbf{Proposition 2.} \emph{If $C$ is a Plott CF, and all CFs $C_x$ are non-empty-valued Plott CFs, then the `integral' CF $D$ is Plottian as well.} \medskip

Proof. Let the element $b$ belong to $B$, but not $D(B)$, and $B'=B-b$. We have to check that $D(B')=D(B)$; this will give the consistency of $D$.

Let's consider two cases: when $b$ is the only element in its fibre, and
when it is not the only one.

If $b$ is the only one, then $\pi (B')=\pi(B)-\pi(b)$. In addition, $\pi(b)$ is not selected from $\pi(B)$. In fact, if $\pi (b)\in C(\pi (B))$, then $b$ as the only element of the fibre $B_{\pi b}$ is selected in this fibre and belongs to  $D(B)$. Due to the consistency of $C$, $C(\pi B')=C(\pi B)$. For points $y$ of $C(\pi B)$, the fibres $B\cap X_y$ and $B'\cap X_y$ are the same, as well as the choices from them. Therefore $D(B')=D(B)$.

If $b$ is not the only one in the fibre $B\cap X_{\pi b}$, then $\pi B'=\pi B$. In all other fibres the choice does not change, so it remains to look at the choice from the fibre over $\pi b$, which we denote as $B_0$ (and the CF
$C_{\pi b}$ as $C_0$). Since $B'_0=B_0-b$, and $b$ is not selected by the function $C_0$, then (due to the consistency of $C_0$) $C_0(B'_0)=C_0(B_0)$. This proves the consistency of $D$.

In the same style, the substitution is checked. Here we assume that $b\in D(B)$ and must check that the remaining elements of $D(B)$ survive when selected from $B'=B-b$. Again, if $\pi B'=\pi B$, then the other fibres do not change and we only need to consider the fibre containing $b$, that is, $Y_{\pi b}$. And in it (due to the substitutability of $C_{\pi b}$), deleting $b$ leads to the fact that everything that was selected continues to be selected.

If $\pi B'$ is less than $\pi B$, then this can happen only if $b$ is a single element of the fibre $B_\pi b$. In this case $\pi B'=\pi B-\pi b$. Also, $\pi b\in C(\pi B)$. From the substitutability of $C$, we get that all other points from $C(\pi B)$ remain in $C(\pi B')$, as well as the choices from the fibres above them. As for the fibre over $\pi b$, this fibre consists only of $b$, which was removed. $\Box$\medskip

\section{Stability problem}

Let's start with a general setup. Let $N$ be the set of \emph{agents} and $E$ be the set of possible \emph{contracts} between agents. For a contract $e$, $P(e)$ denotes the set of \emph{participants} of the contract; symmetrically, for $a\in N$, $E(a)$ is the set of contracts in which the agent $s$ participates. We call this data the \emph{frame} (or skeleton) of the problem. We call a \emph{problem} equipping each agent $a\in N$ with a Plott CF $C_a$ on the set $E(a)$. For the sake of brevity, we say that the agent $a$ is \emph{linear} (respectively, \emph{sequential, cardinally monotonic}, etc.) if the corresponding choice function $C_a$ is linear (respectively, sequential, cardinally monotonic, etc.).\medskip

\textbf{Definition.} The set (or system) of contracts $S\subseteq E$ is called \emph{stable} if

1. For any agent $a\in N$ the set $S(a)=S\cap E(a)$ is acceptable.

2. There are no \emph{blocking} contracts, that is, contracts outside of $S$  which are desirable for all participants of the contract. In other words, if a contract $e$ belongs to $C_a(S(a)\cup e)$ (that is desirable) for every $a\in P(e)$ then $e\in S$.\medskip

We say also that a stable set is a \emph{solution} of the problem.

Our main goal is to disaggregate some agents so that they have become linear. Since a linear agent can conclude only one contract (or remain without contracts) this often greatly simplifies the analysis of the situation.

Let us spend some time to convince the reader that the situation `many-to-one' is really a bit simpler than the general one. Let's assume that our frame is bipartite, that is all agents are divided into `workers' and `firms', and the contracts are paired, between a worker and a firm. We will assume that all workers $w$ are linear; as for firms, we will assume that they have cardinally monotonous Plott CFs. It is well known that in this situation stable systems exist.

We propose here a simple way to build stable systems by induction. More precisely, let us assume that we have a solution $S_0$ in the problem in which a certain worker $w_0$ has been deleted. Then it is possible to canonically rebuild $S_0$ into the solution $S$ of the original problem, already with the worker $w_0$. Below we describe in detail one step of this procedure, inspired by Gale-Shapley algorithm.

Let us assume that we have a problem $I$ with agents $M$ and $W$, a set of
contracts $E$, and choice functions that are linear for workers and cardinally monotonous for firms. Let $I(-w_0)$ denote the problem in which the worker $w_0$ is deleted, and let $S_0$ be some stable set in the problem $I(-w_0)$.

We return the worker $w_0$ together with the bouquet $E(w_0)$ of his contracts. Some of these contracts may be of interest to firms. Say that a contract $e$ from $E(w_0)$ is \emph{interested} for a firm $m$ if $e\in E(m)$ and $e\in C_m(S_0(m)\cup e)$. The set of such contracts is denoted by $D(m)$, and let $D=\cup _m D(m)$.

It may turn out that $D$ is empty (let's call this \emph{situation 0}), that is no firm has shown interest in contracts with $w_0$. Obviously, in this case $S=S_0$ is stable in the  problem $I$.

Otherwise, $D$ is nonempty, and let $d$ be the best (for $w_0$) contract from this set, the contract with some firm $m$.

\emph{Situation 1}. $C_m(S_0(m)\cup d)=S_0(m)\cup d$. Let $S=S_0\cup d$. Again, it is quite obvious that $S$ is stable in the problem $I$. The main thing here is that the set $S(m)=S_0(m)\cup d$ is acceptable for $m$.

\emph{Situation 2.} $C_m(S_0(m)\cup d)$ is strictly less than $S_0(m)\cup d$. Due to the assumption of the cardinal monotonicity of CF $C_m$, the size
$C_m(S_0(m)\cup d)$ cannot be less than the size of $C_m(S_0(m))=S_0(m)$.
Therefore, $C_m(S_0(m)\cup d)$ is obtained from $S_0(m)$ by throwing out exactly one old contract $s$ and the addition of the contract $d$ instead. This rejected contract $s$ was a contract between $m$ and some worker $w_1$ (different, of course, from $w_0$). And if we now delete the worker $w_1$ along with the set $E(w_1)$ of his contracts, that is, consider the problem $I(-w_1)$, then the system $S_1=(S_0-s)\cup d$ will be stable. The fact is that a new desirable contract $d'$ could appear only at the firm $m$. But if $d'$ is desirable in the state $(S_0(m)-s)\cup d$, it is all the more desirable in the state $S_0$ (let's recall the transitivity of the relation $\preceq _m$).\medskip

This completes the description of the first step of the procedure. In the situations 0 and 1, we already got a stable $S$. If the situation 2 is realized, it is necessary to repeat the operation, returning the worker $w_1$, and so on. Sooner or later, this process it will end up in the situation 0 or 1. Indeed, if at the first step we are in the situation 2, the position of the firm $m$ is strictly improving (in the sense of $\preceq_m$), and nothing changes for other firms. In the second step (if we gets in the situation 2 again), there is an improvement for another firm and so on. But firms cannot be improved infinitely, and sooner or later the process gets into situation 0 or 1 and gives a stable state $S$ for the original problem $I$.

Note that for firms, the state $S$ is not worse than the initial one $S_0$. As for the workers, their position fluctuates in different directions during the process, sometimes improving (for $w_0$), sometimes getting worse (for $w_1$). As a kind of curiosity, we note that it may well turn out that in the final stable state $S$ the initiator of the process $w_0$ may turn out to be a loner.

One can see it in the following simple example.

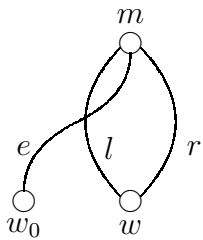
\unitlength=.7mm
\special{em:linewidth 0.4pt}
\linethickness{0.4pt}
\begin{picture}(85.00,50)(-20,0)
\put(70.00,40.00){\circle{4.00}}
\put(70.00,10.00){\circle{4.00}}
\put(50.00,10.00){\circle{4.00}}
\bezier{152}(68.00,11.00)(55.00,25.00)(68.00,39.00)
\bezier{152}(72.00,11.00)(85.00,25.00)(72.00,39.00)
\bezier{76}(50.00,12.00)(50.00,20.00)(60.00,25.00)
\bezier{76}(70.00,38.00)(70.00,30.00)(60.00,25.00)
\put(70.00,45.00){\makebox(0,0)[cc]{$m$}}
\put(70.00,5.00){\makebox(0,0)[cc]{$w$}}
\put(50.00,5.00){\makebox(0,0)[cc]{$w_0$}}
\put(66.00,20.00){\makebox(0,0)[cc]{$l$}}
\put(82.00,20.00){\makebox(0,0)[cc]{$r$}}
\put(50.00,20.00){\makebox(0,0)[cc]{$e$}}
\end{picture}

The firm $m$ ranks the contracts as follows: $l>e>r$, for the worker $w$ on the contrary $r>l$. When $w_0$ is absent, the contract $r$ is stable for the pair $m,w$. But here $w_0$ appears, the firm switches to $e$, abandoning $r$. The worker $w$ temporarily drops out, then returns, the firm offers him the contract $l$, to which $w$ agrees. The contract $e$ is broken, and we get a stable set $\{l\}$. \medskip

In principle, this can be used to improve (for firms) a stable
system of contracts by introducing a `fictitious' worker who disappears at the end of the process. However, I have not yet been able to accurately describe the conditions...

Another benefit of the hypothesis that all workers have linear CFs is that we can compare two systems $S$ and $T$, considering their symmetric difference
$S\triangle T$ and analyzing the resulting rotations. However, so far this is only a direction for further work.

\section{Disaggregation of sequential workers}

Let us assume that there is a distinguished agent (`the worker') 0; and $M=N-\{0\}$. We assume also that $C_0=C_1*C_2$, where $C_1$ and $C_2$ are Plott CFs. Roughly speaking, the agent 0 is an (ordered) pair of subagents 1 and 2. And the idea is to replace the agent 0 by this pair of agents. A new problem $\widetilde I$ arises in which the agents are 1, 2 and agents from $\widetilde M=M$. The new set of contracts $\widetilde E$ is arranged as follows: contracts from $E(0)$ are duplicated, turning into $E(1)$ and $E(2)$ (each is a copy of $E(0$)); the rest (the set $R=E-E(0)$) remains unchanged.

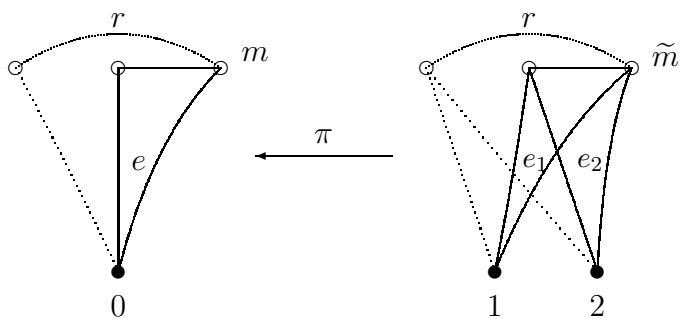
\unitlength=.9mm
\special{em:linewidth 0.4pt}
\linethickness{0.4pt}
\begin{picture}(121.00,55.00)
\put(30.00,40.00){\circle{2}}
\put(45.00,40.00){\circle{2.00}}
\put(60.00,40.00){\circle{2.00}}
\put(90.00,40.00){\circle{2.00}}
\put(105.00,40.00){\circle{2.00}}
\put(120.00,40.00){\circle{2.00}}
\put(115.00,10.00){\circle*{2.00}}
\put(100,10.00){\circle*{2.00}}
\put(45.00,10.00){\circle*{2.00}}
\bezier{70}(30.00,40.00)(45.00,50.00)(60.00,40.00)
\bezier{70}(90.00,40.00)(105.00,50.00)(120.00,40.00)
\put(45.00,5.00){\makebox(0,0)[cc]{0}}
\put(100.00,5.00){\makebox(0,0)[cc]{1}}
\put(115.00,5.00){\makebox(0,0)[cc]{2}}

\put(85,27){\vector(-1,0){20}}
\put(75.00,30.00){\makebox(0,0)[cc]{$\pi$}}

\put(48.00,26.00){\makebox(0,0)[cc]{$e$}}
\put(106.00,26.00){\makebox(0,0)[cc]{\small $e_1$}}
\put(114.00,26.00){\makebox(0,0)[cc]{\small$e_2$}}

\put(45.00,47.00){\makebox(0,0)[cc]{$r$}}
\put(105.00,47.00){\makebox(0,0)[cc]{$r$}}

\put(65.00,42.00){\makebox(0,0)[cc]{$m$}}
\put(125.00,42.00){\makebox(0,0)[cc]{$\widetilde{m}$}}

\bezier{38}(30,40)(38,25)(45,10)

\bezier{38}(90,40)(95,25)(100,10)
\bezier{38}(90,40)(102,25)(115,10)

\bezier{168}(45,40)(47,40)(60,40)
\bezier{168}(45,40)(45,25)(45,10)
\bezier{168}(45,10)(50,30)(60,40)

\bezier{168}(105,40)(112,40)(120,40)
\bezier{168}(105,40)(103,25)(100,10)
\bezier{168}(120,40)(108,30)(100,10)
\bezier{168}(105,40)(110,25)(115,10)
\bezier{168}(120,40)(116,30)(115,10)

\end{picture}

CF for agent 1 is $C_1$, and for agent 2 is $C_2$. For every new firm $\widetilde m$ (an avatar of the old firm $m$) CF $\widetilde C_{\widetilde m}$ is arranged as follows. The set $\widetilde E(\widetilde m)$ is projected into $E(m)=R(m)\coprod E(0,m)$ (where $E(0,m)$ is a set of contracts between 0, $m$ and maybe someone else). Over $R(m)$ we have a bijection. And $E(0,m)$ is duplicated; one double is denoted by $E(1,m)$, the second by $E(2.m)$. CF $\widetilde C_{\widetilde m}$ is constructed as an "integral"\ over CF $C_m$ (see section 3.3), where $\pi: \widetilde{E}(\widetilde{m}) \to E(m)$ is the natural mapping, and in all fibers the level 1 is more preferable than the level 2.

Suppose now that $\widetilde S \subseteq \widetilde E$ is a stable set in the problem $\widetilde I$, and $S=\pi (\widetilde S)$.\medskip

\textbf{Theorem 1.} \emph{The map $\pi :\widetilde S\to S$ is bijective and $S$ is stable set in the problem $I$.}\medskip

The proofs this and the next Theorems  (somewhat tedious, or involved) are postponed to Section 6.

The match $\widetilde  S\mapsto S=\pi (\widetilde S)$ defines the mapping (denoted as $\pi _*$) from ${\bf St}(\widetilde I)$ in ${\bf St}(I)$ (where {\bf St} denotes the set of stable systems).\medskip

\textbf{Theorem 2.} \emph{Mapping $\pi_*:{\bf St}(\widetilde I) \to  {\bf       St}(I)$ is a bijection.}\medskip

Theorem 2 shows that the problems $I$ and $\widetilde I$ are equivalent in the sense that their solutions are canonically identical. This assertion has two natural generalizations.

First, we can assume that CF of agent 0 is decomposed into several CFs. In particular, if CF of the agent 0 is sequential, then this agent can be replaced (as indicated the above method) on an ordered set of linear subagents.

Secondly, we can assume that there are several agents of the type 0. The only thing that needs to be demanded here is that these agents (call them 'workers') are not connected to each other with contracts. Let's say more exactly. Suppose that the set of agents $N$ is divided into two disjoint groups $M$ and $W$, and $E$ is again a set of contracts. We say that the agents from $W$ do not \emph{interact directly} (or are not connected to each other by the contracts) if, for any contract $e\in E$, the intersection $P(e)\cap W$ has  no more than one element. Generalizing Theorem 2, it is easy to get a more general statement:\medskip

\textbf{Theorem 3.} \emph{Suppose that}

1) \emph{$N$ is divided into two parts $M$ and $W$, and the agents from $W$ are not connected to each other;}

2) \emph{for each $w\in W$ the CF $C_w$ is represented as $C_w=C_{w(1)}*...* C_{w(q_w)}$, where all $C_{w(j)}$ are Plott CFs;}

3) \emph{$C_m$ for any $m\in M$ is a Plott CF.}

\emph{Then the initial problem can be replaced by another one, with participants $\widetilde N=\widetilde M\coprod\widetilde W$ and  mappings $\pi : \widetilde N\to N$ and $\pi: \widetilde{E} \to E$, so that}

a) \emph{$\widetilde M$ is bijective $M$, and for $w$ from $W$ the fiber $\pi      ^{-1}(w)$ consist of $q(w)$ elements;}

b) \emph{CF $\widetilde C_{\widetilde w}$ of any agent $\widetilde w=w(j)$ is  $C_{w(j)}$, and}

c) \emph{matching $\widetilde S\mapsto \pi (\widetilde S)$ gives a bijection of the set ${\bf St}(\widetilde I)$ of stable systems in the modified problem and the set ${\bf St}(I)$ of stable systems in the initial problem.}\medskip

Roughly speaking, it is possible to replace each sequential worker $w$ by the set of its linear clones, so that the problems remain equivalent. The proof of  Theorem 3 can be carried out by repeating the reasoning of the Theorems 1 and 2. But it is easier to note that when we split a single agent $w$, the other agents from $W$ don't change at all. Therefore they remain sequential and can be split one by one in turn.\medskip

I would like to finish with one more refinement. The set ${\bf St}(I)$ of stable contract systems is not just a set; different systems can be compared from the point of view of agents. Let's say that, for an agent $m$, a system $S$ is no better than a system $T$ (we write $S\preceq_m T$) if $C_m(S(m)\cup T(m)) \subseteq T(m)$ (see Blair's order from Section 2.1). Similarly, one can talk about systems $\widetilde S$, $\widetilde T$, agents $\widetilde m$ and use the term $\widetilde\preceq_{\widetilde m}$. Well, in the situation of Theorem 2 (or 3), the following holds\medskip

\textbf{Proposition 3.} \emph{If $\widetilde S\ \widetilde \preceq _{\widetilde m}\ \widetilde T$ then $S\preceq _m T$.}\medskip

Proof. By the definition $\widetilde S \ \widetilde \preceq _{\widetilde m} \ \widetilde T$ means that $\widetilde C_{\widetilde m} (\widetilde S(\widetilde m)\cup \widetilde T(\widetilde m))\subseteq \widetilde T(\widetilde m)$. In  force of (\ref{pro}) we have $\pi (\widetilde C_{\widetilde m} (\widetilde S(\widetilde m)\cup \widetilde T(\widetilde m)))=C_m(S(m)\cup T(m))$, and this set is contained in $\pi (\widetilde T(\widetilde m))=T(m)$. That is $S\preceq _m T$.  $\Box$\medskip

If we define the relation  $\widetilde\preceq _{\widetilde M}$ as the intersection of the relations $\widetilde \preceq_{\widetilde m}$ for all $\widetilde m \in \widetilde M$, and similarly shall understand $\preceq_M$, we obtain that
$\widetilde S \ \widetilde \preceq _{\widetilde M} \ \widetilde T$ implies $S \preceq _M T$. This means that the bijection $\pi_*: {\bf St}(\widetilde I)\to {\bf St}(I)$ is a monotone mapping of posets (with partial orders $\widetilde \preceq _{\widetilde M}$ and $\preceq _M$). It is naturally to expect that this is not just a morphism of posets, but an \emph{isomorphism}. However we have only been able to prove this in one important special case.

Namely, let us assume that our framework is bipartite, that is, all agents are divided into two groups: $M$ (`firms') and $W$ (`workers'), and agents within each group are not connected by contracts. This means that all contracts are pair and concluded between agents from different groups. As before, we assume that the workers are sequential. As for CFs of firms, we assume that they are \emph{cardinally monotonous}. With these assumptions, it is true the following refinement of Proposition 3.\medskip

\textbf{Proposition 4.} \emph{Mapping $\pi _*: {\bf St}(\widetilde I) \to {\bf St}(I)$ is an isomorphism of the posets (with orders $\widetilde\preceq =\widetilde \preceq _{\widetilde M}$ and $\preceq =\preceq _M$).}\medskip

Proof. Let $\widetilde{S}$ and $\widetilde{T}$ be two stable sets in the modified problem $\widetilde{I}$. And $S=\pi (\widetilde S)$, $T=\pi (\widetilde T)$ are the corresponding stable sets in $I$. By virtue of Proposition 3, it remains to show that if $S\preceq T$ then $\widetilde S \ \widetilde{\preceq} \ \widetilde T$.

In the modified problem $\widetilde I$ workers are linear and firms are cardinally monotonous (see the remark from Section 2.3). In this situation, the set $({\bf St}(\widetilde I), \widetilde \preceq )$ is not just a poset, but a lattice (\cite{B}). Moreover, the join $\widetilde S\vee \widetilde T$ of two stable systems $\widetilde S$ and $\widetilde T$ is given by the following explicit formula (\cite{A,F}): for any firm $\widetilde m$
                              $$
(\widetilde S\vee \widetilde T)(\widetilde m)=\widetilde C_{\widetilde m} (\widetilde S(\widetilde m)\cup \widetilde T(\widetilde m)).
                                $$
According to (\ref{pro}) we have $\pi ((\widetilde S\vee \widetilde T)(\widetilde m))=C_m(S(m)\cup T(m))=T(m)$, and this holds for any firm $m\in M$. This implies that $\pi_*$ of the stable system $\widetilde S\vee\widetilde T$ is equal to $T$. But by virtue of Theorem 2, this means that $\widetilde S\vee \widetilde T$ is equal to $\widetilde{T}$, that is, $\widetilde S\widetilde \ \preceq \ \widetilde T$.  $\Box$\medskip

Thus, in this particular situation, the solutions of the problems $I$ and $\widetilde I$ are equivalent in a strong sense, not only as sets, but as lattices.

\section{Discussion}

Thus, we have shown that if all workers are sequential, then one can replace them by linear workers. This makes the task of analysis a little easier, because now each worker can  conclude no more than one contract. Workers which left without contracts in one stable system will be the same in any stable system (provided that firms are cardinally monotonous). We cannot delete them, because the removal may lead to the expansion of the set {\bf St} (see the simple example from Section 3). However, even in this case one can delete them, and then see which stable systems disappear after their return. So, in principle, we can assume that the workers conclude exactly one contract.

Of course, it is possible to disaggregate firms instead of workers, if they are sequential. But it is worth emphasizing that the `reduction to linearity'       remains `one-sided': we reduce either workers or firms. To do this simultaneously for both, if possible, it requires some other approach. In \cite{DKar} a simultaneous reduction of all agents was performed, but there the reduced agents had CFs generated by weak orders (rather than linear ones). And although it would seem that the weak orders and linear orders are very close, in the stability problems the difference becomes quite significant. This can be seen on the following example.\medskip

      \textbf{Example.} The frame of the problem looks like this:

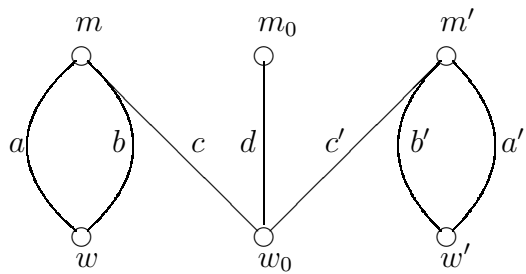
\unitlength=.8mm
\special{em:linewidth 0.4pt}
\linethickness{0.4pt}
\begin{picture}(115.00,55.00)(-10,0)
\put(70.00,40.00){\circle{3.00}}
\put(70.00,10.00){\circle{3.00}}
\put(40.00,40.00){\circle{3.00}}
\put(40.00,10.00){\circle{3.00}}
\bezier{160}(41.00,39.00)(56.00,25.00)(41.00,11.00)
\bezier{172}(39.00,39.00)(23.00,25.00)(39.00,11.00)
\put(100.00,40.00){\circle{3.00}}
\put(100.00,10.00){\circle{3.00}}
\bezier{160}(99.00,11.00)(85.00,25.00)(99.00,39.00)
\bezier{160}(101.00,11.00)(115.00,25.00)(101.00,39.00)

\put(71,11){\line(1,1){28}}
\put(69,11){\line(-1,1){28}}
\put(70,12){\line(0,1){27}}

\put(39,44){$m$}
\put(69,44){$m_0$}
\put(99,44){$m'$}
\put(39,5){$w$}
\put(69,5){$w_0$}
\put(99,5){$w'$}

\put(28,24){$a$}
\put(109,24){$a'$}
\put(45,24){$b$}
\put(94,24){$b'$}
\put(58,24){$c$}
\put(80,24){$c'$}
\put(66,24){$d$}

\end{picture}

Three firms: $m,m'$ and $m_0$, three workers $w,w'$ and $w_0$, and seven contracts $a,b, c, a',b', c'$ and $d$. The preferences of the workers are as follows: $w$ likes $b$ more than $a$, $w'$ likes $b'$ more than $a'$, $w_0$ wants $d$ the least, and $c$ and $c'$ are equivalent for him. For $m$ the contract $a$ is better than $b$ and $c$, which are equivalent for it. Similarly for $m'$: $a'$ is better than $b'$ and $c'$, which are equivalent.

           There are four stable systems in this problem.

           1. The best for firms: $\overline{S}=\{a,a',d\}$.

           2. The best for workers: $\underline{S}=\{b,c,c',b'\}$.

           3. The left: $S_l=\{a,b',c'\}$, and

           4. The right: $S_r=\{b,c,a'\}$.

From the point of view of firms, $\overline{S}$ is the join of $S_l$ and $S_r$. However, if we take the union of $S_l$ and $S_r$, we get a set that does not contain $d$. So their join $\overline{S}$ is not contained in the union. And $\overline{S}$ is not the choice of the firms from $S_l\cup S_r$.\medskip

Blair\cite{B} provides a more complex example where the lattice of stable sets is not distributive.

But all this is a little aside. As for inapplicability of the above construction for reducing both firms and workers to linear ones, we propose  here another simple example. With two agents and two contracts. One of which (say, the left one) is better for both agents. But they are allowed the quota equal to 2. So they use both contracts. If we literally apply the proposed above construction, we get two $M$-agents and two $W$-agents. For all these clones, the left contract is better than the right one, and everyone will choose the left one. So the projection of the solution of modified problem is not a solution in the original one.\medskip

We could proceed as follows. First we `linearize' the workers, and then
(assuming firms are sequential) repeat the construction, but this time in relation to the firms. The failure here is the follows. Even if the firm's CF $C_m$ was sequential, the CF $\widetilde C_{\widetilde m}$ of the modified firm is not longer sequential and even quotable.

However, if the firms were initially linear, then  they remain linear after the proposed modification. Therefore, if the firms were linear, and the workers weew sequential (it is better, of course, to consider the opposite case), then the problem can be reduced to the classical matching case when all agents are linear.

\section{Proofs of Theorems 1 and 2}

\textbf{Proof of Theorem 1.} We need to check the acceptability and non-blocking of the system $S=\pi(\widetilde{S})$. For this we will use the similar properties of the system $\widetilde S$. Let us express these properties in terms of the initial problem $I$.\medskip

\textbf{Lemma 2.} \emph{$\widetilde S$ is acceptable for an agent $\widetilde m$ if and only if $\pi :\widetilde S(\widetilde m) \to S(m)$ bijection and $S$ is acceptable for agent $m$.}\medskip

The Lemma follows from Lemma 1. $\Box$\medskip

Next, what does the acceptability of $\widetilde S$ mean for agents 1 and 2? The answer is given by the following obvious\medskip

\textbf{Lemma 3.} \emph{$\widetilde S$ is acceptable for $i=1,2$ if and only if       $S(i)=C_i(S(i))$. } $\Box$\medskip

To get the acceptability of $S$ for agent 0, we need the following\medskip

\textbf{Lemma 4.} \emph{Let $e\in S(2)$. Then $C_1(S(1)\cup e)=S(1)$.} \medskip

Proof. Due to the consistency of $C_1$, it is sufficient to show that $e \notin C_1(S(1)\cup e)$. Let's suppose that $e\in C_1(S(1)\cup e)$, and denote       by $e_1$ the element $e$, considered as an element of $\widetilde E(1)$. ($e_1$ is $e$ raised on the first floor.) The element $e_1$ does not belong to $\widetilde S(1)$ because it is not in $S(1)$ ($\pi (e_1)=e\in S(1)$). $e_1$ does not belong to $\widetilde S(2)$, because it is located on the ground floor, whereas $\widetilde S(2)$ lies on the second one. Therefore, $e_1$ does not belong to $\widetilde S$.

We claim that $e_1$ blocks $\widetilde S$. To do this, we need to check that $e_1$ is desirable for agent 1, as well as for any agent $\widetilde m$, participating in this contract on the side $\widetilde M$. Desirability for agent 1 is obvious from the fact that $e\in C_1(S(1)\cup e)$. Now let $\widetilde m$ be an agent from $\widetilde M$, involved in this contract. Then $m$ is a participant of the contract $e$. We need to check that $e_1\in \widetilde C_{\widetilde m}(\widetilde S(\widetilde m)\cup e_1)$. Remembering the definition $\widetilde C_{\widetilde m}$, we need to check that $e$ belongs to $C_m(\pi(\widetilde S(\widetilde m)\cup e_1))$. Let $e_2$ be the same contract $e$, but lies on the level 2. Then it belongs to $\widetilde S(2)$ and even more to $\widetilde S$. And to $\widetilde S(\widetilde m)$. Therefore, $\pi (\widetilde S(\widetilde m))$ contains $\pi (e_2)=e$, and $\pi (\widetilde S(\widetilde m)\cup e_1)=\pi (\widetilde S(\widetilde m))\cup \pi (e_1)=\pi (\widetilde S(\widetilde m))$. On the other hand, from the acceptability of $\widetilde S$ for $\widetilde m$ we know (see Lemma 2) that $C_m(\pi(\widetilde S(\widetilde m)))=\pi (\widetilde S(\widetilde m))$ and therefore contains $e$. So, $e_1$ is desirable for all parties of this contract, and besides does not belong to $\widetilde S$. But this contradicts stability $\widetilde S$. This contradiction completes the proof of Lemma 4. $\Box$\medskip

A small comment. Why is $e_1$ desirable for $\widetilde m$? Because its `double' $e_2$ belongs to $\widetilde S(\widetilde m)$, but $e_1$ is `the same contract', only more preferred by the agent 1. Actually, this is just that place where the `preference to negotiate with the clone 1 instead of 2' plays a role.\medskip

\textbf{Corollary.} $S(1)=C_1(S(0))$ \emph{and} $S(2)=C_2(S(0)-C_1(S(0)))$. \medskip

Indeed, if $e$ from $S(0)=S(1)\coprod S(2)$ belongs to $C_1(S(0))$ then (by virtue of the substitutability of $C_1$) $e$ belongs to $C_1(S(1)\cup e)$. By virtue of Lemma 4, in this case $e\in S(1)$. Therefore $C_1(S(0))\subseteq S(1)$. Since $C_1(S(1))=S(1)$, then from the consistency of $C_1$ we get $C_1(S(0))=S(1)$, that is the first equality. The second equality follows from the fact that      $S(0)-C_1(S(0))=S(0)-S(1)=S(2)$ and Lemma 3. $\Box$\medskip

In particular, $C_0(S(0))=S(0)$, which proves the acceptability of $S$ for the agent 0. Thus the acceptability of the system $S$ for all agents is established.\medskip

\emph{Non-blocking.} Now let's check that the set $S$ is non-blocked in the problem $I$. Suppose that some contract $e$ (not belonging to $S$) is desirable for all its participants. It's pretty clear what it could be only the contract with the worker 0. Further, $m$ denotes an arbitrary other       the participant of the contract $e$. Desirability for 0 means that

      1) $e\in C_0(S(0)\cup e)$.

Desirability for $m$ means that

      2) $e\in C_m(S(m)\cup e)$.

Recall that the choice $C_0$ is obtained by first applying $C_1$ to the set      $S(0)\cup e$, and then applying $C_2$ to the remainder.

Let's first consider the case when $e$ is chosen already at the first step,      that is, $e\in C_1(S(0)\cup e)$. Due to the substitutability, $e$ is also chosen from the lesser set of $S(1)\cup e$, $e\in C_1(S(1)\cup e)$. But this means that $e_1$ is desirable for the agent 1. We claim that $e_1$ is also desirable for any agent $\widetilde m$. That is, we must show that $e_1\in \widetilde C_{\widetilde m}(\widetilde S(\widetilde m)\cup e_1)$. Let's recall the definition of $\widetilde C$. We have to descent (using $\pi$) the set $\widetilde S(\widetilde m)\cup e_1$ into $E$, apply $C_m$ to the image, and then "minimally lift"\ the result in $\widetilde S(\widetilde m)\cup e_1$. Projection of $\widetilde S(\widetilde m)\cup e_1$ is equal to $S(m)\cup e$. A choice from this set, according to 2), contains $e$. And the minimal lift of $e$ in $\widetilde S(\widetilde m)\cup e_1$ is exactly equal to $e_1$. So $e_1$ is desirable for any $\widetilde m$. Therefore, $e_1$ blocks $\widetilde S$, what contradicts to the stability of the latter.

This means that $e$ cannot be chosen at in the first step. That is $e\notin C_1(S(0)\cup e)$, from where $C_1(S(0)\cup e)=C_1(S(0))=S(1)$ (the last follows from Corollary). From this we conclude that $e\in C_2(S(0)-S(1)\cup e)=C_2(S(2)\cup e)$. That is, $e_2$ is desirable for the worker 2. Reasoning as above, we see that $e_2$ is desirable for every $\widetilde m$, a participant of the contract $e_2$. Again, this contradicts to stability of $\widetilde S$.

This proves that $S$ is stable. Theorem 1 is proven.  $\Box$\medskip

\textbf{Proof of Theorem 2.} The injectivity is almost obvious. The matter of fact is that $S=\pi (\widetilde S)$ uniquely defines $\widetilde S$. Over $R$ it is a bijection. And for $S(0)$ we have canonical decomposition (due to Corollary) on  $S(1)=C_1(S(0))$ and $S(2)=C_2(S-S(1))$. Which defines $\widetilde S$ over $S(0)$: $S(1)$ goes up to the first floor, and $S(2)$ goes up to the second floor.

\emph{Surjectivity.} The previous remark suggests how to build a mapping  inverse to $\pi_*$. That is, how to define, for each stable set $S$ in the problem $I$, the set $\widetilde S$ in the problem $\widetilde I$. Namely, $S$ consists of $S(0)$ (contracts involving the agent 0) and $R=S-S(0)$. The set $R$ is unambiguously lifted. As for $S(0)$, it decomposes (as above) on $S(1)=C_1(S(0))$ and $S(2)=C_2(S(0)-S(1))$; the first set lifts to the first floor, the second set lifts to the second floor. This gives the set $\widetilde S$ which is projected on $S$. It remains only to show that $\widetilde S$ is stable in the problem $\widetilde I$. That is, to show the acceptability and non-blocking.\medskip

\emph{Acceptability} for $\widetilde m$. By the construction, $\widetilde S$ is bijectively mapped to $S$. And since $S$ is acceptable for all $m$, then by Lemma 2 $\widetilde S$ is acceptable for all $\widetilde m$.

\emph{Acceptability} for agents $i=1, 2$. According to Lemma 3, it is necessary check that $S(i)=C_i(S(i))$ for $i=1,2$. By the definition of substitutability, $S(1)= C_1(S(0))$; in particular, $C_1(S(1))=S(1)$. Finally, $S(2)=C_2(S(0)-S(1))=C_2(S(2))$.

\emph{Non-blocking.} Suppose that some contract $\widetilde e$ from $\widetilde E$ is desirable (under the condition $\widetilde S$) for all its participants. We need to show that $\widetilde e \in \widetilde S$.

Let's start with desirability for $\widetilde m$. Almost obviously from (\ref{pro}) that this means that the contract $e=\pi (\widetilde e)$ is desirable for $m$ (under the condition $S$). If the agents 1 and 2 do not  participate in  $\widetilde e$ (that is, 0 does not participate in the contract $e$), then $e\in S$ and $\widetilde e\in \widetilde S$. Therefore, one can assume that agents 1 or 2 (but not both simultaneously) participate in the contract $\widetilde e$.

Let's start with the case when $\widetilde e$ is desirable for Agent 1. This means that

a) $\widetilde e=e_1$ is located on the ground floor, and

b) $e_1\in \widetilde C_1(\widetilde S(1)\cup e_1)$, which by the definition means $e\in C_1(S(1)\cup e)$.

\noindent By virtue of "path independence" we have
                                                 $$
C_1(S(0)\cup e)=C_1(C_1(S(0))\cup  e)=C_1(S(1)\cup e),
                                                   $$
so $e\in C_1(S(0)\cup e)\subseteq C_0(S(0)\cup e)$. The last means that $e$ is desirable for agent 1. It turns out that $e$ is desirable (under $S$) for all its participants. Due to the stability of $S$, we have $e\in S$ and even $S(1)$. But in this case, $e_1\in\widetilde S$, and we got what we wanted.

Let's now consider the second possibility, when $\widetilde e$ is desirable for worker 2. This means that

a) $\widetilde e=e_2$ is located on the second floor, and

b) $e_2\in\widetilde C_2(\widetilde S(2)\cup e_2)$, that is, $e\in C_2(S(2)\cup e)$.

\noindent  By the definition
                      $$
C_0(S(0)\cup e)=C_1(S(0)\cup e)\cup C_2(S(0)\cup b-C_1(S(0)\cup b).
                        $$
Since $\widetilde e$ is undesirable for the agent 1, $C_1(S(0)\cup e)=S(1)$. Therefore
              $$
C_2(S(0)\cup b-C_1(S(0)\cup b))=C_2(S(0)\cup e-S(1))=C_2(S(2)\cup b)
                $$
and contains $e$. In particular, $e\in C_0(S(0)\cup e)$. Hence the contract $e$ is desirable for all its participants, therefore (due to the stability of $S$) $e\in S$ and even $e\in S(2)$. Therefore, $\widetilde e=e_2$ belongs to $\widetilde S$.

This completes the check of stability of $\widetilde S$ and the proof of Theorem 2. $\Box$

\end{document}